\documentclass[submission,nomarks]{dmtcs-episciences}

\usepackage[utf8]{inputenc}
\usepackage[T1]{fontenc}
\usepackage{amsmath,amssymb}
\usepackage{booktabs,array}
\newcolumntype{L}[1]{>{\raggedright\arraybackslash}p{#1}}
\usepackage[round]{natbib}

\newtheorem{theorem}{Theorem}
\newtheorem{lemma}[theorem]{Lemma}
\newtheorem{corollary}[theorem]{Corollary}
\newtheorem{proposition}[theorem]{Proposition}
\newtheorem{definition}[theorem]{Definition}
\newtheorem{remark}[theorem]{Remark}

\newcommand{\Sym}[1]{S_{#1}}
\newcommand{\id}{\mathrm{id}}
\newcommand{\sgn}{\operatorname{sgn}}
\newcommand{\Heap}{\textup{\textsc{Heap}}}
\newcommand{\ol}[1]{\left[\,#1\,\right]}

\author{Lu\'{\i}s M. S. Russo\thanks{This work was supported by
  national funds through Funda\c{c}\~{a}o para a Ci\^{e}ncia e a
  Tecnologia (FCT), project UIDB/50021/2020, DOI
  \url{https://doi.org/10.54499/UIDB/50021/2020}.}}
\title[Heap-type permutation generation]{Heap-type permutation
  generation with a prescribed final arrangement}
\affiliation{INESC-ID, Instituto Superior T\'{e}cnico,
  Universidade de Lisboa, Portugal}

\begin{document}
\maketitle

\begin{abstract}
Heap's algorithm generates all permutations of an array of n items,
obtaining each one from its predecessor by a single swap. We
determine the net permutation applied to the array by every complete
recursive call: a transposition when the number of items is odd and
at least three, and a cycle through all items when it is even. We
show that the output of each call is a concatenation of right cosets
of the symmetric group on one item fewer, and we prove a general
formula: the order in which the items are parked at the last
position may be chosen freely, every such order corresponds to
exactly one valid schedule, and it determines the net permutation of
the call. Two variants follow. In the cyclic variant a complete
sweep over an even number of items rotates the array by one
position. In the extremal variant, for even n, a sweep started from
the identity ends in the lexicographically largest arrangement that
any single-swap enumeration can leave; by a parity argument this is
the reversal exactly when n is congruent to 2 modulo 4. The
construction of the extremal variant is proved correct when n is
divisible by 4 and verified by computer for all even n up to 10000.
\end{abstract}

\section{Introduction}\label{sec:intro}

\subsection{Motivation}
Listing all permutations of $n$ objects is one of the oldest tasks in
combinatorial computing. It underlies exhaustive search and
brute-force optimisation, the testing of programs and hardware on all
input orders, the verification of conjectures about small permutation
groups, and the enumeration of orbits and of other structures indexed
by permutations. For these uses the cost of \emph{producing} each
permutation matters, which has motivated \emph{minimal-change}
(Gray-code) methods, in which every permutation differs from its
predecessor by a very small modification, so that per-permutation work
is $O(1)$ and the consumer can often update its own data incrementally
instead of recomputing from scratch.

\subsection{Related work}
\paragraph{Algorithms for permutation generation.}
The literature on generating permutations on computers goes back to
the beginning of the electronic-computer era. The early account by
\citet{lehmer1960} of how to teach combinatorial generation to a
computer, the generation by transpositions of \citet{wells1961} and
the procedure \textsc{Permute} of \citet{peck1962} (ACM Algorithm~86)
are among the first published methods. The method of Steinhaus,
Johnson and Trotter, which produces each permutation from the previous
one by an \emph{adjacent} transposition, is described in Algorithm~115
of \citet{trotter1962} and in the note of \citet{johnson1963}; the
name also recalls \citet{steinhaus1964}. The method of
\citet{heap1963}, the subject of this paper, uses a single swap of two
arbitrary positions and a short recursive description. Subsequent work
includes the algorithm of \citet{langdon1967}, the studies of
permutation sequences by \citet{ordsmith1970,ordsmith1971}, the four
algorithms of \citet{ives1976}, and the \emph{loopless} algorithms of
\citet{ehrlich1973}, which perform $O(1)$ work in the worst case per
permutation. Lexicographic enumeration, in which successive
permutations are \emph{not} close in the swap metric, is the basis of
the ``next permutation'' functions of standard libraries (for example
\texttt{std::next\_permutation} in C++; see \citealp{iso14882}) and is
treated in detail by \citet{dijkstra1976} and \citet{knuth4a}. The
survey of \citet{sedgewick1977} compares these and many other methods;
it concluded that Heap's method was, at the time, the most effective
for computer generation of permutations. Textbook treatments are given
by \citet{nijenhuis1978}, \citet{reingold1977}, \citet{kreher1999}
and, most comprehensively, by \citet[Section~7.2.1.2]{knuth4a}.

\paragraph{Gray codes and Cayley graphs.}
A minimal-change listing of permutations is a Hamiltonian path in a
Cayley graph of the symmetric group whose generators are the allowed
``small changes''. Heap's algorithm is such a path for the generating
set of \emph{all} transpositions. \citet{kompel1975} showed that a set
of transpositions whose graph is a tree (for instance a star,
i.e.\ the transpositions $(1\ i)$) already yields a Hamiltonian path,
and the Cayley graph of the star transpositions is the \emph{star
graph} studied as an interconnection network by \citet{akers1989}.
Comprehensive surveys of the area, which also cover many other
combinatorial objects, are due to \citet{savage1997} and
\citet{mutze2023}.

\paragraph{Correctness proofs and the final arrangement.}
Heap's algorithm is usually justified by checking small cases or by an
informal inductive argument; such arguments are available, for
instance, in expository online material such as the Wikipedia article
\citep{wikiheap}, and some of this material also notes properties of
the array left by a complete run (the Wikipedia proof, for example,
speaks of rotations of the array). To our knowledge, however, a
treatment in terms of the \emph{net permutation} of every recursive
call, valid for both parities and expressed with the group structure,
is not available in the literature we examined, although we make no
claim to have searched exhaustively. Such a treatment is more than a
matter of taste: when an enumeration is embedded in a larger
computation (nested loops over several groups, repeated sweeps,
restoring a data structure after a sweep), the \emph{final}
arrangement and its algebraic description are what the surrounding
code needs.

\subsection{Contributions}
\begin{enumerate}
  \item A self-contained primer on permutation notation
    (Section~\ref{sec:notation}): two-line (``matrix'') notation;
    one-line notation (the bottom row only, which is how arrays are
    stored); cycle notation and its \emph{canonical} form, with the
    rules and several examples; how to multiply and invert
    permutations in each notation; and how a swap of two array
    positions acts on them.
  \item A proof by induction, in cycle notation, of the \emph{net
    permutation} of every recursive call and of the structure of the
    output (Theorem~\ref{thm:main}).
  \item A \emph{general net formula} (Theorem~\ref{thm:general}):
    parking the items at the last position in an arbitrary order
    $e_1,\dots,e_{k-1}$ determines a unique valid schedule, and the
    net permutation is $\pi_{k-1}^{\,k}\,(k\ e_1\ \cdots\ e_{k-1})$.
    All previous results are special cases.
  \item Two new variants (Sections~\ref{sec:cyclicvariant}
    and~\ref{sec:lex}): a \emph{cyclic} variant whose net permutation
    for even $k$ is $(1\ 2\ \cdots\ k)$, and a \emph{lexicographically
    extremal} variant for even $n$ that ends in the lexicographically
    largest arrangement attainable by any single-swap enumeration,
    together with a parity bound showing which arrangement that is
    (reversal exactly when $n\equiv2\pmod4$).
  \item A computational verification of all claims that are not proved
    (Section~\ref{sec:verif}): exhaustive runs of the variants for
    small $n$, and a check of the construction of the extremal variant
    for all even $n\le10000$. The verification program and a C
    implementation of the three variants are publicly available
    \citep{permgen}.
\end{enumerate}
We do not claim that the idea of reordering the swaps of Heap's
algorithm is new; what is new, to our knowledge, is the
group-theoretic description of the whole family of Heap-type schedules
and the use of it to prescribe the final arrangement.

\subsection{Organisation}
Section~\ref{sec:notation} is the notation primer.
Section~\ref{sec:alg} defines the algorithm and the admissible
schedules, Section~\ref{sec:results} states the results for them and
Section~\ref{sec:proofs} proves them. Section~\ref{sec:general} proves
the general net formula, and Section~\ref{sec:lex} uses it to
construct the extremal variant. Section~\ref{sec:examples} gives
worked examples, Section~\ref{sec:verif} describes the computational
verification, and Section~\ref{sec:concl} concludes.

\section{Permutation notation}\label{sec:notation}

Let $[n]=\{1,\dots,n\}$ and let $\Sym{n}$ be the set of bijections
$\sigma:[n]\to[n]$; $|\Sym{n}|=n!$. For $k\le n$ we regard
$\Sym{k}\subseteq\Sym{n}$ by letting a permutation of $[k]$ fix
$k+1,\dots,n$. We use two running examples in $\Sym{5}$.

\subsection{Two-line (matrix) notation}\label{sec:twoline}
A permutation is written as a matrix with two rows: the top row lists
$1,2,\dots,n$ and the bottom row lists the corresponding images.
\[
\sigma=\begin{pmatrix}1&2&3&4&5\\3&5&1&2&4\end{pmatrix},\qquad
\tau=\begin{pmatrix}1&2&3&4&5\\2&3&1&5&4\end{pmatrix}.
\]
The first matrix means $\sigma(1)=3$, $\sigma(2)=5$, $\sigma(3)=1$,
$\sigma(4)=2$, $\sigma(5)=4$. Each column is an arrow
$x\mapsto\sigma(x)$, so the columns may be permuted freely without
changing the permutation; for instance
$\sigma=\begin{pmatrix}2&5&1&4&3\\5&4&3&2&1\end{pmatrix}$.

\subsection{One-line notation (bottom row only)}\label{sec:oneline}
Dropping the top row leaves the \emph{one-line notation}
$\ol{\sigma(1)\ \sigma(2)\ \cdots\ \sigma(n)}$:
\[
\sigma=\ol{3\ 5\ 1\ 2\ 4},\qquad \tau=\ol{2\ 3\ 1\ 5\ 4}.
\]
This is exactly how an array of the numbers $1,\dots,n$ is stored in a
program: \texttt{a[i]}$=\sigma(i)$. We use square brackets for
one-line notation to distinguish it from cycle notation.

\subsection{Cycle notation}\label{sec:cycle}
Start at $1$ and follow $x\mapsto\sigma(x)$ until you return to $1$;
this gives one cycle. Repeat from the smallest element not yet used.
For $\sigma$ above: $1\mapsto3\mapsto1$ gives $(1\ 3)$, and
$2\mapsto5\mapsto4\mapsto2$ gives $(2\ 5\ 4)$, so
\[
\sigma=(1\ 3)(2\ 5\ 4),\qquad \tau=(1\ 2\ 3)(4\ 5).
\]
Conventions: a cycle is read cyclically,
$(a_1\ a_2\ \cdots\ a_m)=(a_2\ \cdots\ a_m\ a_1)$; cycles of length
$1$ (fixed points) are omitted, and the identity is written $\id$;
disjoint cycles commute; to find $\sigma(x)$, locate $x$ in its cycle
and take the next entry (cyclically), or $x$ itself if $x$ does not
occur. A cycle of length $m$ is an \emph{$m$-cycle}, a $2$-cycle is a
\emph{transposition}. An $m$-cycle is a product of $m-1$
transpositions, so its sign is $(-1)^{m-1}$; the sign $\sgn$ is
multiplicative.

\subsection{Canonical cycle notation}\label{sec:canon}
The cycle notation of Section~\ref{sec:cycle} is not unique: every
cycle can be rotated, and the cycles can be listed in any order. For
instance $\sigma=\ol{3\ 5\ 1\ 2\ 4}$ can be written $(1\ 3)(2\ 5\ 4)$,
$(3\ 1)(2\ 5\ 4)$, $(2\ 5\ 4)(1\ 3)$ or $(5\ 4\ 2)(3\ 1)$. A
\emph{canonical} form fixes these choices, so that every permutation
has exactly one representation; it makes permutations comparable by
inspection, and it is the form in which the net permutations of this
paper are tabulated. We use the following rules.
\begin{description}
  \item[(C1)] Rotate every cycle so that it \emph{begins with its
    smallest element}.
  \item[(C2)] List the cycles in \emph{increasing order of their first
    (smallest) elements}.
  \item[(C3)] Omit cycles of length $1$ (fixed points); the identity
    is written $\id$. (If a ``complete'' form is wanted, the fixed
    points are written as $1$-cycles in their places according to
    (C2).)
\end{description}
The representation is unique because the smallest element of a cycle
is unique and the cycles are disjoint, so their smallest elements are
distinct and (C2) determines their order. The procedure of
Section~\ref{sec:cycle} (start from the smallest element not yet used,
follow images until the cycle closes) produces (C1)--(C3)
automatically. To \emph{test} whether a written product of disjoint
cycles is canonical, check that every cycle starts with its minimum
and that these minima increase from left to right.

\paragraph{Examples.}
\begin{itemize}
  \item $\sigma=\ol{3\ 5\ 1\ 2\ 4}$ has the canonical form
    $(1\ 3)(2\ 5\ 4)$; the three other writings above violate (C1) or
    (C2). For $\tau=\ol{2\ 3\ 1\ 5\ 4}$ the canonical form is
    $(1\ 2\ 3)(4\ 5)$.
  \item The products computed in Section~\ref{sec:mult} are canonical
    as written: $\sigma\tau=(1\ 5\ 2)$ and $\tau\sigma=(2\ 4\ 3)$.
    Note that $(2\ 4\ 3)$ is canonical (it starts with its minimum
    $2$) although it is not increasing.
  \item The net permutation of classic Heap on $k=6$ items is
    $(1\ 6\ 5\ 2\ 3\ 4)$, canonical as written. The expression
    $(b_1\ \cdots\ b_{k-1}\ k)$ of Theorem~\ref{thm:main} is a cycle
    that ends with $k$; its canonical form is obtained by rotating it
    to start with $1$.
  \item The net permutation of the extremal variant
    (Section~\ref{sec:lex}) for $n=8$ is $(1\ 7\ 2\ 8)(3\ 6)(4\ 5)$ in
    canonical form, while $(8\ 1\ 7\ 2)(6\ 3)(5\ 4)$ is the same
    permutation written non-canonically.
\end{itemize}

\paragraph{Other conventions.}
``Canonical cycle notation'' is not used in the same sense by all
authors, so the convention should always be stated. The fundamental
transformation of \citet{foata1968} \citep[see
also][\S1.3]{stanleyEC1} writes every cycle \emph{starting with its
largest} element and lists the cycles in \emph{increasing order of
their largest elements}; erasing the parentheses then gives a one-line
permutation. Other authors keep the smallest element first but list
the cycles in \emph{decreasing} order of these elements. For $\sigma$
above, the three conventions give
\[
\begin{array}{ll}
(1\ 3)(2\ 5\ 4) & \text{(this paper)},\\
(3\ 1)(5\ 4\ 2) & \text{(largest first, increasing)},\\
(2\ 5\ 4)(1\ 3) & \text{(smallest first, decreasing)},
\end{array}
\]
and the corresponding one-line words after erasing parentheses are
$\ol{3\ 1\ 5\ 4\ 2}$ and $\ol{2\ 5\ 4\ 1\ 3}$ for the last two. To
pass from our form to the largest-first form, rotate each cycle to
start with its maximum and sort the cycles by their first entries. Our
convention is the natural one for algorithms that discover cycles by
scanning $1,2,\dots,n$ and for results stated, like ours, with a
distinguished point (the point~$1$, or the point~$k$).

\subsection{Multiplication (composition)}\label{sec:mult}
We compose right to left: $(\sigma\tau)(x)=\sigma(\tau(x))$,
i.e.\ $\tau$ acts first. (Some books compose left to right; with our
convention products read like function composition.) Products are
associative but not commutative. For our running example:
\[
\sigma\tau=\ol{5\ 1\ 3\ 4\ 2}=(1\ 5\ 2),\qquad
\tau\sigma=\ol{1\ 4\ 2\ 3\ 5}=(2\ 4\ 3).
\]
We show how each is obtained in each notation.

\paragraph{One-line notation.} Use the entries of $\tau$ as indices
into $\sigma$: $(\sigma\tau)[i]=\sigma[\tau[i]]$. Here
$\sigma[\tau[1]]=\sigma[2]=5$, $\sigma[\tau[2]]=\sigma[3]=1$,
$\sigma[\tau[3]]=\sigma[1]=3$, $\sigma[\tau[4]]=\sigma[5]=4$,
$\sigma[\tau[5]]=\sigma[4]=2$, giving $\ol{5\ 1\ 3\ 4\ 2}$. (In C:
\texttt{c[i] = s[t[i]]}.)

\paragraph{Two-line notation.} Stack three rows: $x$, then $\tau(x)$
(the bottom row of $\tau$), then $\sigma(\tau(x))$, obtained by
looking up the middle entry in the top row of $\sigma$ and reading
what is below it. The product is the first and last rows:
\[
\begin{pmatrix}1&2&3&4&5\\2&3&1&5&4\\5&1&3&4&2\end{pmatrix}
\ \Longrightarrow\ \sigma\tau=
\begin{pmatrix}1&2&3&4&5\\5&1&3&4&2\end{pmatrix}.
\]
Equivalently, reorder the columns of $\sigma$ so that its top row
equals the bottom row of $\tau$ and glue.

\paragraph{Cycle notation.} Trace every point through the factors,
from the rightmost factor to the leftmost. For
$\sigma\tau=(1\ 3)(2\ 5\ 4)\cdot(1\ 2\ 3)(4\ 5)$: the point $1$ goes
$1\to2$ in $(1\ 2\ 3)$, $2\to5$ in $(2\ 5\ 4)$, and $5$ is untouched
by $(1\ 3)$; so $1\mapsto5$. Next, $5\to4$ in $(4\ 5)$, $4\to2$ in
$(2\ 5\ 4)$, $2$ is untouched by $(1\ 3)$; so $5\mapsto2$. Then
$2\to3\to3\to1$, so $2\mapsto1$, which closes the cycle $(1\ 5\ 2)$.
The points $3$ and $4$ return to themselves ($3\to1\to1\to3$ and
$4\to5\to4\to4$).

\subsection{Inversion}\label{sec:inv}
The inverse $\sigma^{-1}$ undoes $\sigma$;
$(\sigma\tau)^{-1}=\tau^{-1}\sigma^{-1}$.
\begin{itemize}
  \item \emph{Two-line:} swap the two rows, then sort the columns by
    the new top row. For $\sigma$: the swapped matrix has columns
    $3\!\to\!1$, $5\!\to\!2$, $1\!\to\!3$, $2\!\to\!4$, $4\!\to\!5$;
    sorted,
    $\sigma^{-1}=\begin{pmatrix}1&2&3&4&5\\3&4&1&5&2\end{pmatrix}$.
  \item \emph{One-line:} for every $i$ write $i$ at index $\sigma[i]$,
    i.e.\ $\sigma^{-1}[\sigma[i]]=i$ (in C: \texttt{inv[s[i]] = i}).
    For $\sigma=\ol{3\ 5\ 1\ 2\ 4}$: $1$ goes to index $3$, $2$ to
    index $5$, $3$ to index $1$, $4$ to index $2$, $5$ to index $4$,
    giving $\ol{3\ 4\ 1\ 5\ 2}$.
  \item \emph{Cycle:} write every cycle backwards:
    $\sigma^{-1}=(1\ 3)(4\ 5\ 2)=(1\ 3)(2\ 4\ 5)$.
\end{itemize}
Check:
$(\sigma\tau)^{-1}=(1\ 5\ 2)^{-1}=(2\ 5\ 1)=(1\ 2\ 5)
=\ol{2\ 5\ 3\ 4\ 1}$,
and $\tau^{-1}\sigma^{-1}$ computed by the one-line rule gives the
same array.

\begin{table}[ht]
\centering
\renewcommand{\arraystretch}{1.25}
\begin{tabular}{@{}L{2.6cm}L{4.0cm}L{2.7cm}L{2.9cm}@{}}
\toprule
 & Two-line & One-line & Cycles\\
\midrule
Product $\sigma\tau$ ($\tau$ first) & stack $\tau$ over $\sigma$ (rows
$x$, $\tau(x)$, $\sigma(\tau(x))$), keep first and last rows &
\texttt{c[i]=s[t[i]]} & trace each point right to left\\
Inverse $\sigma^{-1}$ & swap rows, re-sort columns &
\texttt{inv[s[i]]=i} & write each cycle backwards\\
Left product $(i\ j)\sigma$ & swap the \emph{values} $i,j$ in the
bottom row & swap the \emph{values} $i,j$ & see
Section~\ref{sec:transp}\\
Right product $\sigma(i\ j)$ & swap the bottom entries of columns $i$
and $j$ (top row unchanged) & swap the \emph{entries} at positions
$i,j$ & see Section~\ref{sec:transp}\\
\bottomrule
\end{tabular}
\caption{How to multiply and invert in each
notation.}\label{tab:cheat}
\end{table}

\subsection{Multiplying by a transposition}\label{sec:transp}
Because Heap's algorithm only ever swaps two positions, products with
a transposition $(i\ j)$ are the important case. In one-line notation:
\begin{itemize}
  \item $(i\ j)\,\sigma$ \emph{(left product)} swaps the \emph{values}
    $i$ and $j$, wherever they occur:
    $(2\ 4)\ol{3\ 5\ 1\ 2\ 4}=\ol{3\ 5\ 1\ 4\ 2}$.
  \item $\sigma\,(i\ j)$ \emph{(right product)} swaps the
    \emph{entries at positions} $i$ and $j$:
    $\ol{3\ 5\ 1\ 2\ 4}(2\ 4)=\ol{3\ 2\ 1\ 5\ 4}$.
\end{itemize}
Indeed $((i\ j)\sigma)(x)=(i\ j)(\sigma(x))$ changes the values,
whereas $(\sigma(i\ j))(x)=\sigma((i\ j)(x))$ reads $\sigma$ at
swapped indices. In cycle notation, multiplying by a transposition
either splits a cycle (if $i$ and $j$ lie in the same cycle) or merges
two cycles (if they lie in different cycles).

\subsection{Permutations acting on arrays; states}\label{sec:action}
An array $a=\ol{a_1\ \cdots\ a_n}$ has positions $1,\dots,n$. A
permutation $\sigma$ \emph{acts} on arrays by moving the item at
position $p$ to position $\sigma(p)$: the new array $b$ satisfies
$b_{\sigma(p)}=a_p$. Performing $\sigma_1$ first and $\sigma_2$ second
has total effect $\sigma_2\sigma_1$. Exchanging the items at positions
$i$ and $j$ is the action of the transposition $(i\ j)$, which is its
own inverse.

\begin{lemma}[Arrays are inverses]\label{lem:arrayinv}
Applying $\sigma$ to the identity array $\ol{1\ 2\ \cdots\ n}$
produces the array
\[
\ol{\sigma^{-1}(1)\ \cdots\ \sigma^{-1}(n)},
\]
that is, the one-line notation of $\sigma^{-1}$.
\end{lemma}
\begin{proof}
The new array satisfies $b_{\sigma(p)}=p$ for all $p$, that is,
$b_q=\sigma^{-1}(q)$.
\end{proof}

At any moment of an algorithm that starts from an array $x$ and only
moves items, the current array is obtained from $x$ by a unique
permutation $\omega$, which we call the \emph{state} (relative to
$x$). The initial state is $\id$. Swapping positions $i$ and $j$
replaces the state $\omega$ by $(i\ j)\,\omega$ (a left product). If
$x$ is the identity array, the current array is $\ol{\omega^{-1}}$;
consistently, swapping the entries at positions $i,j$ of
$\ol{\omega^{-1}}$ gives
$\omega^{-1}(i\ j)=\bigl((i\ j)\omega\bigr)^{-1}$ by
Section~\ref{sec:transp}.

\begin{table}[ht]
\centering
\renewcommand{\arraystretch}{1.15}
\begin{tabular}{@{}cccccc@{}}
\toprule
step $m$ & swap & array $\ol{\omega_m^{-1}}$ & state $\omega_m$
(cycles) & state $\omega_m$ (one-line) &
$\omega_m=(\text{swap})\,\omega_{m-1}$\\
\midrule
0 & --- & $\ol{1\ 2\ 3}$ & $\id$ & $\ol{1\ 2\ 3}$ & ---\\
1 & $(1\ 2)$ & $\ol{2\ 1\ 3}$ & $(1\ 2)$ & $\ol{2\ 1\ 3}$ &
$(1\ 2)\,\id$\\
2 & $(1\ 3)$ & $\ol{3\ 1\ 2}$ & $(1\ 2\ 3)$ & $\ol{2\ 3\ 1}$ &
$(1\ 3)(1\ 2)$\\
3 & $(1\ 2)$ & $\ol{1\ 3\ 2}$ & $(2\ 3)$ & $\ol{1\ 3\ 2}$ &
$(1\ 2)(1\ 2\ 3)$\\
4 & $(1\ 3)$ & $\ol{2\ 3\ 1}$ & $(1\ 3\ 2)$ & $\ol{3\ 1\ 2}$ &
$(1\ 3)(2\ 3)$\\
5 & $(1\ 2)$ & $\ol{3\ 2\ 1}$ & $(1\ 3)$ & $\ol{3\ 2\ 1}$ &
$(1\ 2)(1\ 3\ 2)$\\
\bottomrule
\end{tabular}
\caption{The six arrays produced by Heap's algorithm for $n=3$ and the
corresponding states. Each new state is the previous one multiplied
\emph{on the left} by the swap. In one-line notation this swaps two
\emph{values} of the state, and equivalently two \emph{entries} of the
array. The last state is $(1\ 3)$, the net permutation of
$\Heap(3)$.}\label{tab:n3}
\end{table}

\paragraph{Example (Table~\ref{tab:n3}).} Going from step $1$ to step
$2$: the state is multiplied on the left by $(1\ 3)$. In cycle
notation, $(1\ 3)(1\ 2)$: $1\to2\to2$, $2\to1\to3$, $3\to3\to1$, so
the product is $(1\ 2\ 3)$. In one-line notation, the state
$\ol{2\ 1\ 3}$ has its \emph{values} $1$ and $3$ exchanged, giving
$\ol{2\ 3\ 1}$. On the array, the entries at positions $1$ and $3$ of
$\ol{2\ 1\ 3}$ are exchanged, giving
$\ol{3\ 1\ 2}=\ol{(2\ 3\ 1)^{-1}}$.

\section{Heap's algorithm and its variants}\label{sec:alg}

\subsection{The algorithm}
Positions are numbered from $1$. The procedure $\Heap(k)$ permutes the
first $k$ positions of a global array $a$:
\begin{enumerate}
  \item if $k=1$: output $a$ and return;
  \item call $\Heap(k-1)$;
  \item for $i=1,2,\dots,k-1$: swap the items at positions $s_k(i)$
    and $k$, then call $\Heap(k-1)$.
\end{enumerate}
The numbers $s_k(i)\in\{1,\dots,k-1\}$ form the \emph{schedule}; the
call $\Heap(n)$ on an array of length $n$ is the whole algorithm. In
\emph{classic} Heap's algorithm
\[
s_k(i)=\begin{cases}1&k\text{ odd},\\ i&k\text{ even}.\end{cases}
\]
Between two consecutive outputs exactly one swap is performed: after
an output the recursion returns to the first loop that still has
rounds left, performs that round's swap, and the next output follows
without any further swap, because every call $\Heap(k')$ begins with a
call to $\Heap(k'-1)$ and so on down to $\Heap(1)$. We will prove that
all $n!$ outputs are distinct, so the algorithm performs $n!-1$ swaps.

\subsection{Admissible schedules}
Write $\iota_2=\id$ and $\iota_k=(1\ k{-}1)$ for even $k\ge4$.

\begin{definition}[Admissible schedule]\label{def:adm}
A schedule $s=(s_k)_{k\ge2}$ is \emph{admissible} if
\begin{description}
  \item[(A1)] $s_k(i)=1$ for all $i$, for every odd $k\ge3$;
  \item[(A2)] for every even $k\ge2$, the sequence
    $b_k=(b_k(1),\dots,b_k(k-1))$ defined by
  \[
  b_k(i)=\begin{cases}\iota_k\bigl(s_k(i)\bigr)& i\text{
  odd},\\ s_k(i)& i\text{ even},\end{cases}
  \]
  is a permutation of $\{1,\dots,k-1\}$.
\end{description}
\end{definition}
Condition (A2) says that $s_k$ is obtained from some permutation $b_k$
by applying the transposition $\iota_k$ at the odd indices; so there
are exactly $(k-1)!$ admissible choices of $s_k$ for each even $k$.

\begin{definition}[Cyclic variant]
The \emph{cyclic variant} uses the odd-$k$ rule of the classic
algorithm and, for even $k\ge4$, the schedule
\[
s_k=(k{-}1,\ 2,\ 3,\ \dots,\ k{-}2,\ 1),
\]
that is, $s_k(1)=k-1$, $s_k(k-1)=1$ and $s_k(i)=i$ otherwise;
$s_2=(1)$. Thus it performs the swaps
$(k{-}1\ k),(2\ k),(3\ k),\dots,(k{-}2\ k),(1\ k)$: the first and last
partners of the classic schedule are interchanged.
\end{definition}

\begin{lemma}\label{lem:bothadm}
The classic and the cyclic schedules are admissible. For even $k\ge4$
the classic schedule has $b_k=(k{-}1,2,3,\dots,k{-}2,1)$ and the
cyclic one has $b_k=(1,2,3,\dots,k{-}1)$; for $k=2$ both have
$b_2=(1)$.
\end{lemma}
\begin{proof}
(A1) is built in. Classic: $b_k(i)=\iota_k(i)$ for odd $i$, and
$b_k(i)=i$ for even $i$; since $\iota_k$ swaps only $1$ and $k-1$
(both odd), $b_k=(k{-}1,2,\dots,k{-}2,1)$. Cyclic:
$b_k(1)=\iota_k(k-1)=1$, $b_k(k-1)=\iota_k(1)=k-1$, and $b_k(i)=i$ for
the remaining indices (for odd $i$ with $3\le i\le k-3$ because
$\iota_k$ fixes $i$, for even $i$ by definition), so
$b_k=(1,2,\dots,k-1)$.
\end{proof}

\subsection{Net permutation and output states}
Fix an admissible schedule. Start $\Heap(k)$ on an array $x$ and let
$\omega$ denote the state relative to $x$ (Section~\ref{sec:action});
the swaps only involve positions $\le k$, so $\omega\in\Sym{k}$, and
they do not depend on the contents of the array.

\begin{definition}
The \emph{net permutation} $\pi_k\in\Sym{k}$ is the state at the end
of a complete call $\Heap(k)$. The \emph{output states}
$\sigma_k^{(1)},\dots,\sigma_k^{(k!)}$ are the states at which the
successive outputs occur (so $\sigma_k^{(1)}=\id$).
\end{definition}

Since the effect of a call does not depend on the contents, a call
started in state $\theta$ (relative to some earlier array) ends in
state $\pi_k\theta$ and its outputs occur at the states
$\sigma_k^{(j)}\theta$.

\section{Results}\label{sec:results}

\begin{theorem}[Structure of Heap-type algorithms]\label{thm:main}
Let $s$ be an admissible schedule. For every $k\ge1$:
\begin{enumerate}[(i)]
  \item \emph{Net permutation.} $\pi_1=\id$; for odd $k\ge3$,
    $\pi_k=(1\ k)$; for even $k\ge2$,
  \[
  \pi_k=\bigl(b_k(1)\ \ b_k(2)\ \ \cdots\ \ b_k(k{-}1)\ \ k\bigr),
  \]
  a $k$-cycle.
  \item \emph{Output.} The output states
    $\sigma_k^{(1)},\dots,\sigma_k^{(k!)}$ are pairwise distinct,
    hence a listing of all of $\Sym{k}$; consecutive output states
    differ by a transposition ($\sigma_k^{(j+1)}(\sigma_k^{(j)})^{-1}$
    is a transposition); and the last one is $\sigma_k^{(k!)}=\pi_k$.
  \item \emph{Block structure.} For $k\ge2$ the output of $\Heap(k)$
    consists of $k$ consecutive blocks, the $m$-th ($m=0,\dots,k-1$)
    being the right coset $\Sym{k-1}\theta_m$ listed in the order of
    $\Heap(k-1)$, where $\theta_m$ is the state at the start of the
    block. During block $m$ the item that started at position $e_m$
    stays at position $k$, where
  \[
  e_0=k,\qquad e_m=\begin{cases}b_k(m)& k\text{
  even},\\ \pi_{k-1}^{-m}(1)& k\text{ odd},\end{cases}\qquad 1\le m\le
  k-1,
  \]
  and $e_0,\dots,e_{k-1}$ are a permutation of $\{1,\dots,k\}$.
\end{enumerate}
\end{theorem}

\begin{corollary}[Classic Heap's algorithm]\label{cor:classic}
$\pi_2=(1\ 2)$; $\pi_k=(1\ k)$ for odd $k\ge3$; and for even $k\ge4$
\[
\pi_k=(1\ \ k\ \ k{-}1\ \ 2\ \ 3\ \ \cdots\ \ k{-}2).
\]
\end{corollary}

\begin{corollary}[Cyclic variant]\label{cor:cyclic}
For the cyclic variant, $\pi_k=(1\ k)$ for odd $k\ge3$ and
$\pi_k=(1\ 2\ 3\ \cdots\ k)$ for even $k\ge2$. Consequently, for even
$k$ a complete call $\Heap(k)$ rotates the array by one position to
the right:
$\ol{a_1\ a_2\ \cdots\ a_k}\mapsto\ol{a_k\ a_1\ \cdots\ a_{k-1}}$.
\end{corollary}

\begin{corollary}[Every $k$-cycle occurs]\label{cor:allcycles}
Let $k\ge2$ be even. The map $s_k\mapsto\pi_k$ is a bijection from the
$(k-1)!$ admissible choices of $s_k$ onto the set of all $k$-cycles of
$\Sym{k}$.
\end{corollary}

\begin{proposition}[Parity obstruction]\label{prop:parity}
Let $k\ge2$. In any procedure that outputs all $k!$ permutations with
consecutive outputs related by a single swap, the last state is an odd
permutation. In particular it is never a $k$-cycle for odd $k$ (a
$k$-cycle has sign $(-1)^{k-1}=+1$). For even $k$ the net permutation
of $\Heap(k)$ is always odd, and the $k$-cycles are odd.
\end{proposition}

\begin{table}[ht]
\centering
\renewcommand{\arraystretch}{1.2}
\begin{tabular}{@{}clclc@{}}
\toprule
$k$ & classic net $\pi_k$ & final array & cyclic net $\pi_k$ & final
array\\
\midrule
2 & $(1\ 2)$ & $\ol{2\ 1}$ & same & same\\
3 & $(1\ 3)$ & $\ol{3\ 2\ 1}$ & same & same\\
4 & $(1\ 4\ 3\ 2)$ & $\ol{2\ 3\ 4\ 1}$ & $(1\ 2\ 3\ 4)$ &
$\ol{4\ 1\ 2\ 3}$\\
5 & $(1\ 5)$ & $\ol{5\ 2\ 3\ 4\ 1}$ & same & same\\
6 & $(1\ 6\ 5\ 2\ 3\ 4)$ & $\ol{4\ 5\ 2\ 3\ 6\ 1}$ &
$(1\ 2\ 3\ 4\ 5\ 6)$ & $\ol{6\ 1\ 2\ 3\ 4\ 5}$\\
7 & $(1\ 7)$ & $\ol{7\ 2\ 3\ 4\ 5\ 6\ 1}$ & same & same\\
8 & $(1\ 8\ 7\ 2\ 3\ 4\ 5\ 6)$ & $\ol{6\ 7\ 2\ 3\ 4\ 5\ 8\ 1}$ &
$(1\ 2\ \cdots\ 8)$ & $\ol{8\ 1\ 2\ 3\ 4\ 5\ 6\ 7}$\\
\bottomrule
\end{tabular}
\caption{Net permutations (Corollaries~\ref{cor:classic}
and~\ref{cor:cyclic}) and the array obtained when a complete call
$\Heap(k)$ is run on $\ol{1\ 2\ \cdots\ k}$. By
Lemma~\ref{lem:arrayinv} the final array is the one-line notation of
$\pi_k^{-1}$. These entries were also confirmed by running the
algorithm.}\label{tab:nets}
\end{table}

\section{Proofs}\label{sec:proofs}

\subsection{Three lemmas on products}

\begin{lemma}[Conjugation]\label{lem:conj}
For $\gamma\in\Sym{n}$ and a transposition $(i\ j)$:
$\ \gamma\,(i\ j)\,\gamma^{-1}=(\gamma(i)\ \ \gamma(j))$.
\end{lemma}
\begin{proof}
Let $x\in[n]$ and $y=\gamma^{-1}(x)$. If $y=i$, then $x=\gamma(i)$ and
$\gamma(i\ j)\gamma^{-1}(x)=\gamma(j)$; if $y=j$ similarly
$\gamma(j)\mapsto\gamma(i)$. If $y\notin\{i,j\}$, then
$x\notin\{\gamma(i),\gamma(j)\}$ and the left-hand side fixes $x$.
\end{proof}

\begin{lemma}[Star product]\label{lem:star}
Let $r\ge1$, let $a_1,\dots,a_r,z$ be distinct points. Then
\[
(a_1\ z)(a_2\ z)\cdots(a_r\ z)=(z\ \ a_r\ \ a_{r-1}\ \ \cdots\ \ a_1),
\]
the cycle $z\mapsto a_r\mapsto a_{r-1}\mapsto\cdots\mapsto a_1\mapsto
z$.
\end{lemma}
\begin{proof}
Induction on $r$. For $r=1$ both sides equal $(a_1\ z)$. Let $r\ge1$
and assume the claim for $r$; we prove it for $r+1$ points
$a_1,\dots,a_{r+1}$. By the induction hypothesis applied to
$a_2,\dots,a_{r+1}$,
\[
(a_1\ z)(a_2\ z)\cdots(a_{r+1}\ z)=(a_1\ z)\,\gamma,\qquad
\gamma=(z\ \ a_{r+1}\ \ a_r\ \ \cdots\ \ a_2).
\]
Apply $\gamma$ first and then $(a_1\ z)$. The point $z$ goes to
$a_{r+1}$, which $(a_1\ z)$ fixes. A point $a_j$ with $3\le j\le r+1$
goes to $a_{j-1}$, which $(a_1\ z)$ fixes since $j-1\ge2$. The point
$a_2$ goes to $z$ and then to $a_1$. The point $a_1$ is fixed by
$\gamma$ and then goes to $z$. So $z\mapsto a_{r+1}\mapsto
a_r\mapsto\cdots\mapsto a_2\mapsto a_1\mapsto z$, as claimed.
\end{proof}

\begin{lemma}[Inserting a point into a cycle]\label{lem:insert}
Let $k\ge3$, let $C$ be a $(k-1)$-cycle on $\{1,\dots,k-1\}$ (fixing
$k$) and $\rho=(1\ k)$. Then $C\rho$ is a $k$-cycle, namely
$C\rho=(1\ \ k\ \ C(1)\ \ C^2(1)\ \cdots\ C^{k-2}(1))$. In particular
$(C\rho)^k=\id$.
\end{lemma}
\begin{proof}
Compute $C\rho$ on each point (apply $\rho$ first): $1\mapsto k\mapsto
k$; $k\mapsto1\mapsto C(1)$; and $x\mapsto C(x)$ for $x\notin\{1,k\}$.
Starting from $1$ we get $1\mapsto k\mapsto C(1)$; from then on, as
long as the current point is neither $1$ nor $k$, the map simply
applies $C$, so we continue $C(1)\mapsto C^2(1)\mapsto\cdots$. Since
$C$ is a $(k-1)$-cycle, the points $C(1),\dots,C^{k-2}(1)$ are exactly
the $k-2$ elements of $\{2,\dots,k-1\}$, and
$C(C^{k-2}(1))=C^{k-1}(1)=1$. So the orbit of $1$ contains all $k$
points and $C\rho$ is a single $k$-cycle.
\end{proof}

\subsection{Unwinding the recursion}

\begin{lemma}[Unwinding]\label{lem:unwind}
Let $k\ge2$, $C=\pi_{k-1}$, and $\tau_i=(s_k(i)\ \ k)$ for $1\le i\le
k-1$. Define $\theta_0=\id$ and $\theta_m=\tau_m\,C\,\theta_{m-1}$ for
$1\le m\le k-1$. Then in a call $\Heap(k)$ started in state $\id$:
\begin{enumerate}[(a)]
  \item the $(m{+}1)$-st of the $k$ calls to $\Heap(k-1)$ (``block
    $m$'') starts in state $\theta_m$, its outputs occur at the states
    $\sigma_{k-1}^{(j)}\theta_m$ ($1\le j\le(k-1)!$), and it ends in
    state $C\theta_m$;
  \item $\theta_m=\tau_m C\,\tau_{m-1}C\cdots\tau_1 C$, and
  \begin{equation}\label{eq:unwind}
  \pi_k=C\,\theta_{k-1}=C\,\tau_{k-1}\,C\,\tau_{k-2}\,C\cdots
  C\,\tau_1\,C;
  \end{equation}
  \item during block $m$ the item at position $k$ is the one that
    started at position $e_m:=\theta_m^{-1}(k)$.
\end{enumerate}
\end{lemma}
\begin{proof}
(a) Block $0$ starts in state $\id=\theta_0$. A block that starts in
state $\theta_m$ ends in state $C\theta_m$ by the remark at the end of
Section~\ref{sec:alg}; then the swap $\tau_{m+1}$ is performed and the
next block starts in $\tau_{m+1}C\theta_m=\theta_{m+1}$. The statement
about the output states is the same remark. (b) is immediate from the
recursion; $\pi_k$ is the state at the end of block $k-1$, namely
$C\theta_{k-1}$. (c) Every state in block $m$ has the form
$\sigma\theta_m$ with $\sigma\in\Sym{k-1}$, so $\sigma$ fixes $k$. The
item now at position $k$ is the item $p$ with $\sigma\theta_m(p)=k$,
i.e.\ $\theta_m(p)=k$, i.e.\ $p=\theta_m^{-1}(k)$.
\end{proof}

\begin{lemma}[Gluing]\label{lem:glue}
Let $k\ge2$ and suppose \textup{(ii)} of Theorem~\ref{thm:main} holds
for $k-1$. If $e_0,\dots,e_{k-1}$ (from Lemma~\ref{lem:unwind}) are
pairwise distinct, then \textup{(ii)} holds for $k$ and the output of
$\Heap(k)$ is the concatenation of the blocks of (iii).
\end{lemma}
\begin{proof}
\emph{Distinctness.} By Lemma~\ref{lem:unwind}(a), block $m$ consists
of the $(k-1)!$ states $\sigma\theta_m$, which are pairwise distinct
because the $\sigma=\sigma_{k-1}^{(j)}$ are. If
$\omega=\sigma\theta_m$ is in block $m$ then
$\omega^{-1}(k)=\theta_m^{-1}\sigma^{-1}(k)=\theta_m^{-1}(k)=e_m$. As
the $e_m$ are distinct, no state lies in two blocks. So there are
$k\cdot(k-1)!=k!$ distinct output states, i.e.\ all of $\Sym{k}$.

\emph{Consecutive states.} Inside block $m$, consecutive states
$\sigma^{(j)}\theta_m$ and $\sigma^{(j+1)}\theta_m$ have quotient
$\sigma^{(j+1)}(\sigma^{(j)})^{-1}$, a transposition by hypothesis.
Between blocks $m$ and $m+1$, the last state of block $m$ is
$\sigma_{k-1}^{((k-1)!)}\theta_m=C\theta_m$ (hypothesis (ii): the last
state is $\pi_{k-1}=C$), and the first state of block $m+1$ is
$\sigma_{k-1}^{(1)}\theta_{m+1}=\theta_{m+1}=\tau_{m+1}(C\theta_m)$,
so the quotient is the transposition $\tau_{m+1}$.

\emph{Last state.} The last output of $\Heap(k)$ is the last output of
its last block, at state $C\theta_{k-1}=\pi_k$ by~\eqref{eq:unwind}.
\end{proof}

\subsection{Proof of Theorem~\ref{thm:main}}

We prove (i), (ii), (iii) together by induction on $k$. The case $k=1$
is trivial: $\Heap(1)$ outputs once, at the identity, so $\pi_1=\id$,
(ii) holds, and (iii) is empty. Let $k\ge2$ and assume the theorem for
$k-1$.

\subsubsection*{Case $k$ odd, $k\ge3$}
Then $k-1$ is even and
$C:=\pi_{k-1}=\bigl(b_{k-1}(1)\ \cdots\ b_{k-1}(k-2)\ \ k{-}1\bigr)$
is a $(k-1)$-cycle on $\{1,\dots,k-1\}$, because $b_{k-1}$ is a
permutation of $\{1,\dots,k-2\}$ (for $k=3$ this reads $C=(1\ 2)$). By
(A1) all swaps equal $\rho=(1\ k)$, so $\tau_i=\rho$ for all $i$.

\emph{Proof of (i).} By~\eqref{eq:unwind}, $\pi_k=(C\rho)^{k-1}C$. By
Lemma~\ref{lem:insert}, $(C\rho)^k=\id$, hence
$(C\rho)^{k-1}=(C\rho)^{-1}=\rho^{-1}C^{-1}=\rho\,C^{-1}$ and
\[
\pi_k=\rho\,C^{-1}C=\rho=(1\ k).
\]

\emph{Proof of (iii).} Here $\theta_m=(\rho C)^m$, so
$\theta_m^{-1}=(C^{-1}\rho)^m$ and $e_m=\psi^m(k)$ with
$\psi:=C^{-1}\rho$. We have $\psi(k)=C^{-1}(1)$, $\psi(1)=k$ and
$\psi(x)=C^{-1}(x)$ otherwise. We claim $e_m=C^{-m}(1)$ for $1\le m\le
k-1$. For $m=1$ this is $\psi(k)=C^{-1}(1)$. If $e_m=C^{-m}(1)$ with
$1\le m\le k-2$, then $e_m\notin\{1,k\}$ (as $C$ has order $k-1$ and
does not fix $1$, we have $C^{-m}(1)\ne1$ for $1\le m\le k-2$), so
$e_{m+1}=\psi(e_m)=C^{-1}(e_m)=C^{-(m+1)}(1)$. The points $C^{-m}(1)$,
$1\le m\le k-1$, are the $k-1$ distinct elements of $\{1,\dots,k-1\}$
(the orbit of $1$ under the $(k-1)$-cycle $C^{-1}$), and $e_0=k$; so
$e_0,\dots,e_{k-1}$ are pairwise distinct.

\emph{Proof of (ii).} Apply Lemma~\ref{lem:glue}.

\subsubsection*{Case $k$ even, $k\ge2$}
By the induction hypothesis (i), $C:=\pi_{k-1}$ equals $\id$ if $k=2$
(since $\pi_1=\id$) and equals $(1\ k{-}1)$ if $k\ge4$ (since $k-1$ is
odd and $\ge3$). In both cases $C=\iota_k$: $C$ is an involution,
$C^{-1}=C$, and $C(k)=k$. Write $s_i=s_k(i)$, $b_i=b_k(i)$ and
$\tau_i=(s_i\ k)$. By Lemma~\ref{lem:conj},
\begin{equation}\label{eq:conjC}
C\,\tau_i\,C=C\,(s_i\ k)\,C^{-1}=\bigl(C(s_i)\ \ k\bigr).
\end{equation}

\emph{Proof of (i).} The word in~\eqref{eq:unwind} is
$C\tau_{k-1}C\tau_{k-2}\cdots C\tau_1C$: it has $k$ letters $C$ and
$k-1$ letters $\tau$, alternating, beginning and ending with $C$.
Since $k-1$ is odd, the indices $i=k-1,k-2,\dots,1$ alternate odd,
even, \dots, odd. Cut the word into consecutive pieces: $C\tau_iC$ for
every odd $i$, and the single letter $\tau_i$ for every even $i$:
\[
\pi_k=(C\tau_{k-1}C)\,\tau_{k-2}\,(C\tau_{k-3}C)\,\tau_{k-4}
\cdots\tau_2\,(C\tau_1C).
\]
By~\eqref{eq:conjC}, the piece belonging to $i$ is the transposition
$f_i=(b_i\ \ k)$ with $b_i=C(s_i)$ for odd $i$ and $b_i=s_i$ for even
$i$; this is exactly the sequence $b_k$ of (A2), a permutation of
$\{1,\dots,k-1\}$. Thus
\[
\pi_k=(b_{k-1}\ k)(b_{k-2}\ k)\cdots(b_1\ k).
\]
Apply Lemma~\ref{lem:star} with $r=k-1$, $z=k$ and $a_j=b_{k-j}$ (so
$a_r=b_1$, $a_{r-1}=b_2,\dots,a_1=b_{k-1}$): the product equals
$(k\ \ b_1\ \ b_2\ \cdots\ b_{k-1})=(b_1\ \ \cdots\ \ b_{k-1}\ \ k)$,
as claimed.

\emph{Proof of (iii).} Since all letters are involutions,
$\theta_m^{-1}=C\tau_1C\tau_2\cdots C\tau_m$ by
Lemma~\ref{lem:unwind}(b) (reverse the order of the factors and invert
each). Cut this word as before into $f_1=C\tau_1C$, $f_2=\tau_2$,
$f_3=C\tau_3C,\dots$: if $m$ is even we get
$\theta_m^{-1}=f_1f_2\cdots f_{m-1}\,g_m$ with $g_m=\tau_m$; if $m$ is
odd we get $\theta_m^{-1}=f_1\cdots f_{m-1}\,g_m$ with a last piece
$g_m=C\tau_m$ (there is no closing $C$). Here $f_j=(b_j\ k)$ as
in~\eqref{eq:conjC}. Now $g_m(k)=s_m=b_m$ if $m$ is even and
$g_m(k)=C(s_m)=b_m$ if $m$ is odd. Each $f_j$ with $j<m$ moves only
$b_j$ and $k$, and $b_m\notin\{b_1,\dots,b_{m-1},k\}$ by (A2), so
$f_1\cdots f_{m-1}$ fixes $b_m$. Hence $e_m=\theta_m^{-1}(k)=b_m$ for
$1\le m\le k-1$, and $e_0=k$. These are pairwise distinct.

\emph{Proof of (ii).} Apply Lemma~\ref{lem:glue}. \qed

\subsection{Proofs of the corollaries and of
Proposition~\ref{prop:parity}}

\begin{proof}[of Corollaries~\ref{cor:classic} and~\ref{cor:cyclic}]
By Lemma~\ref{lem:bothadm} both schedules are admissible, so
Theorem~\ref{thm:main}(i) applies. Classic, even $k\ge4$:
$b_k=(k{-}1,2,3,\dots,k{-}2,1)$, so $\pi_k=(k{-}1\ \ 2\ \ 3\cdots
k{-}2\ \ 1\ \ k)$, which after rotating to start at $1$ is
$(1\ k\ k{-}1\ 2\ 3\ \cdots\ k{-}2)$. Cyclic: $b_k=(1,2,\dots,k-1)$,
so $\pi_k=(1\ 2\ \cdots\ k)$. This permutation sends position $p$ to
$p+1$ for $p<k$ and $k$ to $1$; the item at position $p$ moves to
$p+1$ and the item at position $k$ moves to position $1$, which is the
rotation claimed. For $k=2$ both give $(1\ 2)$.
\end{proof}

\begin{proof}[of Corollary~\ref{cor:allcycles}]
By Theorem~\ref{thm:main}(i), $\pi_k=(b_k(1)\ \cdots\ b_k(k{-}1)\ k)$.
A $k$-cycle of $\Sym{k}$ can be written in exactly one way as a
sequence that ends with $k$ (read the cycle cyclically starting just
after $k$), so the $(k-1)!$ permutations $b_k$ correspond bijectively
to the $(k-1)!$ $k$-cycles. The map $b_k\mapsto s_k$ given by
$s_k(i)=\iota_k(b_k(i))$ for odd $i$ and $s_k(i)=b_k(i)$ for even $i$
is a bijection (as $\iota_k$ is an involution) onto the admissible
choices of $s_k$.
\end{proof}

\begin{proof}[of Proposition~\ref{prop:parity}]
The last state is a product of $k!-1$ transpositions, one for each of
the $k!-1$ swaps (consecutive outputs differ by one swap and the last
output occurs after the last swap). Each transposition has sign $-1$
and $\sgn$ is multiplicative, so the sign is $(-1)^{k!-1}=-1$ because
$k!$ is even for $k\ge2$. A $k$-cycle has sign $(-1)^{k-1}$, which is
$+1$ for odd $k$.
\end{proof}

\begin{remark}
The proof shows where the odd/even dichotomy comes from. For odd $k$
the same swap $(1\ k)$ is used in every round and the argument
collapses to a power of a single $k$-cycle $C\rho$: a $k$-cycle has
order $k$, so the $k-1$ swaps telescope to $\pi_k=(1\ k)$. For even
$k$ the inner net permutation is an involution, and conjugating by it
turns the $k-1$ rounds into a product of star transpositions, which
Lemma~\ref{lem:star} evaluates as a single $k$-cycle. The choice of
the schedule affects only the order of the points
$b_k(1),\dots,b_k(k-1)$ in that cycle.
\end{remark}

\section{A general net formula: parking orders}\label{sec:general}

Section~\ref{sec:alg} restricted the schedules by the admissibility
conditions (A1) and (A2). In this section we drop all restrictions and
describe \emph{every} schedule at one level by the order in which
items are ``parked'' at the last position. The result
(Theorem~\ref{thm:general}) contains Theorem~\ref{thm:main} as a
special case and is the tool used in Section~\ref{sec:lex}.

\subsection{Setting}
Fix $k\ge2$ and a permutation $C\in\Sym{k-1}\subseteq\Sym{k}$ (so
$C(k)=k$); in applications $C=\pi_{k-1}$, but nothing below uses this.
Let $s=(s(1),\dots,s(k-1))$ be a schedule for level $k$ (each
$s(i)\in\{1,\dots,k-1\}$), $\tau_i=(s(i)\ \ k)$, and define as in
Lemma~\ref{lem:unwind}
\begin{align*}
&\theta_0=\id,\qquad \theta_m=\tau_m\,C\,\theta_{m-1}\ (1\le m\le
k-1),\\
&e_m=\theta_m^{-1}(k)\ (0\le m\le k-1),\qquad
\pi(C,s)=C\,\theta_{k-1}.
\end{align*}
Thus $e_0=k$ and $e_m$ is the item (identified by its starting
position) parked at position $k$ during block $m$. We call $s$
\emph{valid for $C$} if $e_0,\dots,e_{k-1}$ are pairwise distinct. By
the proof of Lemma~\ref{lem:glue} this is exactly the condition that
the $k$ blocks lie in different right cosets of $\Sym{k-1}$; if
$C=\pi_{k-1}$ and the schedules below level $k$ are valid, then $s$ is
valid for $C$ if and only if $\Heap(k)$ outputs every permutation of
$\Sym{k}$ exactly once. (If two parked items coincide,
$\theta_m\theta_{m'}^{-1}$ fixes $k$ and two blocks lie in the same
coset, so outputs repeat.)

\begin{lemma}[Path product]\label{lem:path}
For pairwise distinct points $e_0,\dots,e_m$:
$\ (e_0\ e_1)(e_1\ e_2)\cdots(e_{m-1}\ e_m)=(e_0\ e_1\ \cdots\ e_m)$.
\end{lemma}
\begin{proof}
Induction on $m$; $m=1$ is trivial. Write the product for $m$ as
$P\,(e_{m-1}\ e_m)$ with $P=(e_0\ \cdots\ e_{m-1})$ by hypothesis, and
apply $(e_{m-1}\ e_m)$ first. Then $e_m\mapsto e_{m-1}\mapsto e_0$;
$e_{m-1}\mapsto e_m\mapsto e_m$ ($P$ fixes $e_m$); and $e_i\mapsto
e_i\mapsto e_{i+1}$ for $i<m-1$. So $e_0\mapsto
e_1\mapsto\cdots\mapsto e_m\mapsto e_0$.
\end{proof}

\begin{theorem}[General net formula]\label{thm:general}
Let $k\ge2$, $C\in\Sym{k-1}$ and let $s$ be any schedule for level
$k$.
\begin{enumerate}[(a)]
  \item $s(m)=C\bigl(\theta_{m-1}(e_m)\bigr)$ for all $m$, and
    $\theta_m=C^m\,(e_0\ e_1)(e_1\ e_2)\cdots(e_{m-1}\ e_m)$.
  \item If $s$ is valid for $C$, then
    $\theta_m=C^m\,(k\ \ e_1\ \cdots\ e_m)$ and
  \[
  \pi(C,s)\;=\;C^{\,k}\,\bigl(k\ \ e_1\ \ e_2\ \cdots\ e_{k-1}\bigr).
  \]
  \item The map $s\mapsto(e_1,\dots,e_{k-1})$ is a bijection from the
    schedules valid for $C$ onto the orderings of $\{1,\dots,k-1\}$.
    In particular there are exactly $(k-1)!$ valid schedules.
  \item Consequently the possible net permutations of level $k$ are
    exactly the elements $C^{\,k}Z$ with $Z$ a $k$-cycle of $\Sym{k}$,
    and a target $T\in\Sym{k}$ is the net permutation of some valid
    schedule if and only if $C^{-k}T$ is a $k$-cycle.
\end{enumerate}
\end{theorem}

\begin{proof}
(a) Since $\theta_m^{-1}=\theta_{m-1}^{-1}C^{-1}\tau_m$ and
$\tau_m(k)=s(m)$, we get
$e_m=\theta_{m-1}^{-1}\bigl(C^{-1}(s(m))\bigr)$, i.e.\ $s(m)=C(u)$
with $u:=\theta_{m-1}(e_m)$. Because $C(k)=k$ and $s(m)\neq k$ we have
$u\ne k$, hence $e_m\ne e_{m-1}$ (recall $\theta_{m-1}(e_{m-1})=k$).
By Lemma~\ref{lem:conj}, $\tau_m=(C(u)\ \ k)=C\,(u\ \ k)\,C^{-1}$, so
\[
\theta_m=\tau_mC\theta_{m-1}=C\,(u\ \ k)\,\theta_{m-1}.
\]
Again by Lemma~\ref{lem:conj}, and since $\theta_{m-1}(e_m)=u$ and
$\theta_{m-1}(e_{m-1})=k$, we have
$(u\ \ k)=\theta_{m-1}(e_{m-1}\ \ e_m)\theta_{m-1}^{-1}$, so
$(u\ \ k)\theta_{m-1}=\theta_{m-1}(e_{m-1}\ \ e_m)$ and
$\theta_m=C\,\theta_{m-1}\,(e_{m-1}\ \ e_m)$. Induction from
$\theta_0=\id$ gives the stated product.

(b) If $e_0,\dots,e_{k-1}$ are distinct, Lemma~\ref{lem:path} turns
the product into the cycle $(e_0\ \cdots\ e_m)=(k\ e_1\ \cdots\ e_m)$,
and $\pi(C,s)=C\theta_{k-1}=C^{k}(k\ e_1\ \cdots\ e_{k-1})$.

(c) Injectivity: by (a), $s(m)$ is determined by $e_m$ and
$\theta_{m-1}$, which in turn is determined by $e_1,\dots,e_{m-1}$.
Surjectivity: given an ordering $(e_1,\dots,e_{k-1})$ of
$\{1,\dots,k-1\}$, define $s(m)=C(\theta_{m-1}(e_m))$ and $\theta_m$
recursively. Suppose the parked items of the first $m-1$ blocks are
the prescribed ones; then $\theta_{m-1}(e_m)\ne k$, because
$\theta_{m-1}(e_{m-1})=k$ and $e_m\neq e_{m-1}$, so
$s(m)=C(\theta_{m-1}(e_m))\in\{1,\dots,k-1\}$ (as $C$ fixes $k$ and
permutes $\{1,\dots,k-1\}$), and by the computation in (a) the item
parked in block $m$ is $e_m$. Hence the schedule exists, and its
parked items are the prescribed, pairwise distinct ones.

(d) By (b) and (c), the nets are $C^kZ$ with $Z=(k\ e_1\cdots
e_{k-1})$, and as $(e_1,\dots,e_{k-1})$ runs over all orderings, $Z$
runs over all $k$-cycles of $\Sym{k}$ (a $k$-cycle has a unique cyclic
representation starting at $k$). So $T=C^kZ$ for a $k$-cycle $Z$ if
and only if $C^{-k}T$ is a $k$-cycle.
\end{proof}

\begin{corollary}[Counting]\label{cor:count}
There are exactly $1!\cdot2!\cdots(n-1)!$ valid Heap-type schedule
families for $n$ items, i.e.\ choices of one valid schedule at each
level $2,\dots,n$.
\end{corollary}
\begin{proof}
By Theorem~\ref{thm:general}(c) there are $(k-1)!$ valid schedules at
level $k$ for each possible lower net $C$, independently of $C$;
multiply over $k=2,\dots,n$.
\end{proof}

\subsection{The cyclic variant in terms of parking
orders}\label{sec:cyclicvariant}
\begin{corollary}[Parking orders of the classic and the cyclic
variant]\label{cor:orders}
\begin{enumerate}[(i)]
  \item For even $k$, $C=\pi_{k-1}$ is an involution fixing $k$, so
    $C^k=\id$ and $\pi_k=(k\ e_1\ \cdots\ e_{k-1})$. The classic
    variant parks $e=(k{-}1,2,3,\dots,k{-}2,1)$, the cyclic variant
    parks $e=(1,2,\dots,k{-}1)$ (increasing order).
  \item For odd $k\ge3$, the classic schedule $s\equiv1$ parks
    $e_m=C^{-m}(1)$ (proof of Theorem~\ref{thm:main}(iii)). If the
    levels below are cyclic, $C=(1\ 2\ \cdots\ k{-}1)$ and so
    $e=(k{-}1,k{-}2,\dots,1)$ (decreasing order), and indeed
    $\pi_k=C\,(k\ e_1\cdots e_{k-1})=(1\ k)$.
\end{enumerate}
Hence the cyclic variant is the Heap-type algorithm that parks the
items in \emph{increasing} order at even levels and in
\emph{decreasing} order at odd levels.
\end{corollary}
\begin{proof}
(i) is Theorem~\ref{thm:main}(i),(iii) rewritten (the parking orders
are $e_m=b_k(m)$, Lemma~\ref{lem:bothadm}). (ii) The first claim is in
the proof of Theorem~\ref{thm:main}; $C^{-m}(1)=k-m$ for
$C=(1\ 2\ \cdots\ k{-}1)$ and $1\le m\le k-1$. For the product, let
$Z=(k\ k{-}1\ k{-}2\ \cdots\ 1)$, i.e.\ $Z(k)=k-1$, $Z(x)=x-1$ for
$2\le x\le k-1$, and $Z(1)=k$. Since $C$ has order $k-1$,
Theorem~\ref{thm:general}(b) gives $\pi_k=C^kZ=CZ$. Now
$CZ(k)=C(k-1)=1$, $CZ(1)=C(k)=k$, and for $2\le x\le k-1$,
$CZ(x)=C(x-1)=x$. So $CZ=(1\ k)$.
\end{proof}

\section{A variant with lexicographically extremal final
arrangement}\label{sec:lex}

The cyclic variant of Corollary~\ref{cor:cyclic} ends in a rotation of
the array. It is natural to ask for an enumeration whose final
arrangement, started from $\ol{1\ 2\ \cdots\ n}$, is the
\emph{lexicographically largest} one, so that a complete sweep lists
the arrangements ``from the smallest to the largest''. We first show
which arrangements are possible at all, then give a construction for
even $n$.

\subsection{Which final arrangements are possible}
Recall (Lemma~\ref{lem:arrayinv}) that if the state at the end is
$\omega$ then the final array is $\ol{\omega^{-1}}$, a permutation of
the same sign as $\omega$.

\begin{proposition}[Lexicographically largest attainable
arrangement]\label{prop:lexbound}
Let $n\ge2$. For every procedure that lists all arrangements of an
$n$-array, starting from $\ol{1\ 2\ \cdots\ n}$, with consecutive
arrangements differing by one swap, the final array is an odd
permutation. The lexicographically largest odd permutation is
\[
a^*_n=\begin{cases}\ol{n\ \ n{-}1\ \cdots\ 2\ \ 1}\ \text{(the
reversal)} & n\equiv2,3\pmod 4,\\[2pt]
\ol{n\ \ n{-}1\ \cdots\ 3\ \ 1\ \ 2}& n\equiv0,1\pmod4.\end{cases}
\]
\end{proposition}
\begin{proof}
By Proposition~\ref{prop:parity} the final state is odd, hence so is
the final array. The reversal is a product of $\lfloor n/2\rfloor$
disjoint transpositions $(i\ \ n{+}1{-}i)$, so it is odd if and only
if $\lfloor n/2\rfloor$ is odd, i.e.\ $n\equiv2,3\pmod4$; then it is
the lexicographically largest permutation of all and in particular the
largest odd one. Otherwise the reversal is even; the permutation that
immediately precedes it in lexicographic order is obtained by
exchanging its last two entries, i.e.\ $\ol{n\ \cdots\ 3\ 1\ 2}$,
which has the opposite sign and is therefore odd, hence is the largest
odd permutation.
\end{proof}

Thus for even $n$ the true reversal can be the final arrangement
\emph{only when $n\equiv2\pmod 4$}; for $n\equiv0\pmod4$
(e.g.\ $n=4,8,12$) the best possible is $\ol{n\ \cdots\ 3\ 1\ 2}$. By
Theorem~\ref{thm:main} no \emph{admissible} schedule can reach the
reversal for $n\ge4$, since the nets of admissible schedules are
$n$-cycles while the reversal is a product of $\lfloor n/2\rfloor\ge2$
disjoint transpositions; the construction below therefore uses a level
that is not admissible. For even $n$ write $T_n=(a^*_n)^{-1}$ for the
target \emph{state}: $T_n$ is the reversal $v\mapsto n{+}1{-}v$ if
$n\equiv2\pmod4$ and
\[
T_n:\ 1\mapsto n{-}1,\quad 2\mapsto n,\quad v\mapsto n{+}1{-}v\ (3\le
v\le n)\qquad\text{if }n\equiv0\pmod4 .
\]

\subsection{The construction}
\begin{definition}[Extremal variant, even $n\ge4$]\label{def:lexmax}
Levels $2,\dots,n-2$ use the cyclic schedule. Level $n-1$ (odd) parks
the items $1,\dots,n-2$ in the order
\[
E_n=\begin{cases}(1,2,\dots,n-2)&
n\equiv0\pmod4,\\
(1,\dots,m,\ n{-}2,\ n{-}3,\ \dots,\ m{+}1),\ m=\tfrac{n-2}{2}&
n\equiv2\pmod4.\end{cases}
\]
Level $n$ parks the items in the order $e_1,\dots,e_{n-1}$ given by
the $n$-cycle
\[
Z_n=(n\ \ e_1\ \cdots\ e_{n-1}):=\pi_{n-1}^{-n}\,T_n ,
\]
provided this permutation is an $n$-cycle. For $n=2$ the variant is
the cyclic (= classic) one. The schedules $s_{n-1}$ and $s_n$ are
obtained from these orders by Theorem~\ref{thm:general}(a):
$s(m)=C(\theta_{m-1}(e_m))$.
\end{definition}
Thus the extremal variant coincides with the cyclic variant except at
the two topmost levels, and level $n$ is defined by one explicit
product of permutations; there is no need to search.

\begin{theorem}[Extremal variant]\label{thm:lexmax}
Let $n\ge2$ be even.
\begin{enumerate}[(a)]
  \item If $Z_n$ is an $n$-cycle, then the extremal variant lists all
    $n!$ permutations with one swap between consecutive arrangements
    and its final array, started from $\ol{1\ 2\ \cdots\ n}$, is
    $a^*_n$: the lexicographically largest arrangement that any
    single-swap enumeration can leave
    (Proposition~\ref{prop:lexbound}).
  \item $Z_n$ is an $n$-cycle whenever $n\equiv0\pmod4$.
  \item $Z_n$ is an $n$-cycle for every even $n$ with $6\le n\le
    10000$ (verified by computer, including $n\equiv2\pmod4$).
\end{enumerate}
\end{theorem}

\begin{proof}[of (a)]
By Corollary~\ref{cor:cyclic} and Theorem~\ref{thm:main}, levels
$2,\dots,n-2$ are valid and $\pi_{n-2}=A:=(1\ 2\ \cdots\ n{-}2)$. By
Theorem~\ref{thm:general}(c) every parking order at level $n-1$ is
valid for $C=A$, so by Lemma~\ref{lem:glue} levels up to $n-1$ are
valid. If $Z_n$ is an $n$-cycle, Theorem~\ref{thm:general}(d) gives a
valid schedule at level $n$ with net $\pi_{n-1}^{\,n}Z_n=T_n$, and
Lemma~\ref{lem:glue} shows that $\Heap(n)$ enumerates $\Sym{n}$ with
one swap between consecutive arrangements. The final state is $T_n$
and by Lemma~\ref{lem:arrayinv} the final array is $T_n^{-1}=a^*_n$.
\end{proof}

\begin{lemma}[Net of level $n-1$ for $n\equiv0\pmod 4$]\label{lem:U}
Let $n\ge4$, $n\equiv0\pmod4$, and let levels $\le n-2$ be cyclic and
level $n-1$ use the order $E_n=(1,\dots,n-2)$. Then
\[
U:=\pi_{n-1}=(1\ 3\ 5\ \cdots\ n{-}3)\,(2\ 4\ \cdots\ n{-}2\ \ n{-}1).
\]
\end{lemma}
\begin{proof}
With $k=n-1$, $C=A=(1\ 2\ \cdots\ n{-}2)$ (of order $n-2$) and
$Z'=(n{-}1\ \ 1\ \ 2\ \cdots\ n{-}2)$, Theorem~\ref{thm:general}(b)
gives $U=A^{n-1}Z'=AZ'$. Now $Z'(x)=x+1$ for $x\le n-2$ and
$Z'(n-1)=1$, and $A(y)=y+1$ for $y\le n-3$, $A(n-2)=1$, $A(n-1)=n-1$.
Hence $U(x)=x+2$ for $1\le x\le n-4$; $U(n-3)=A(n-2)=1$;
$U(n-2)=A(n-1)=n-1$; $U(n-1)=A(1)=2$. The odd points $1,3,\dots,n-3$
form one cycle, and $2,4,\dots,n-2,n-1$ another.
\end{proof}

\begin{proof}[of (b)]
For $n=4$ we have $U=(2\ 3)$ and $U^4=\id$, so $Z_4=T_4=(1\ 3\ 2\ 4)$
is a $4$-cycle. Let $n\ge8$. Write $V=(1\ 3\ \cdots\ n{-}3)$, of
length $L=\frac{n-2}{2}$, and $W=(2\ 4\ \cdots\ n{-}2\ \ n{-}1)$, of
length $\frac n2$. Since $n=2\cdot\frac n2=2L+2$, $W^n=\id$ and
$V^n=V^2$, so $Z_n=U^{-n}T_n=V^{-2}T_n$. The permutation $V^{-2}$
moves only odd points: $y\mapsto y-4$ for odd $y\ge5$, $3\mapsto n-3$,
$1\mapsto n-5$. With the explicit $T_n$ above, $Z:=Z_n$ is therefore
given by
\begin{align*}
&Z(1)=n{-}1,\quad Z(n{-}1)=2,\quad Z(2)=n,\quad Z(n)=n{-}5,\\
&Z(x)=n+1-x\quad\text{for odd }x,\ 3\le x\le n-3,\\
&Z(x)=n-3-x\quad\text{for even }x,\ 4\le x\le n-4,\qquad
Z(n{-}2)=n{-}3.
\end{align*}
(For odd $x$ in the range, $T_n(x)=n+1-x$ is even and fixed by
$V^{-2}$; for even $x\le n-4$, $T_n(x)$ is odd $\ge5$ and $V^{-2}$
subtracts $4$; for $x=n-2$, $T_n(x)=3$.) Put $J=\frac{n-4}{4}$ and,
for suitable $j$,
\[
a_j=n-1-4j,\quad b_j=4j+2\quad(1\le j\le J);\qquad c_j=n-3-4j,\quad
d_j=4j+4\quad(0\le j\le J-1).
\]
Then $Z(n)=a_1$; $Z(a_j)=b_j$ for $1\le j\le J$ (as $3\le a_j\le
n-3$); $Z(b_j)=a_{j+1}$ for $j\le J-1$ (as $4\le b_j\le n-4$);
$Z(b_J)=Z(n-2)=n-3=c_0$; $Z(c_j)=d_j$ for $0\le j\le J-1$ (as $5\le
c_j\le n-3$); $Z(d_j)=c_{j+1}$ for $j\le J-2$; and
$Z(d_{J-1})=Z(n-4)=1$. Hence the orbit of $1$ is
\[
1\to n{-}1\to2\to n\to a_1\to b_1\to\cdots\to a_J\to b_J\to c_0\to
d_0\to\cdots\to c_{J-1}\to d_{J-1}\to1 .
\]
It consists of $4+4J=n$ points, and they are pairwise distinct: modulo
$4$ the points $1,c_0,\dots,c_{J-1}$ are $\equiv1$ (the values
$1,5,9,\dots,n{-}3$), the points $n{-}1,a_1,\dots,a_J$ are $\equiv3$
(the values $3,7,\dots,n{-}1$), the points $2,b_1,\dots,b_J$ are
$\equiv2$ (the values $2,6,\dots,n{-}2$) and the points
$n,d_0,\dots,d_{J-1}$ are $\equiv0$ (the values $4,8,\dots,n$). So the
orbit is all of $\{1,\dots,n\}$ and $Z_n$ is an $n$-cycle.
\end{proof}

\paragraph{The case $n\equiv2\pmod4$.}
Here the same reduction applies with a different explicit $U$: for
$m=\frac{n-2}{2}$ one finds (and we have verified by computer for
$n=6,10,\dots,2000$) that
\[
\pi_{n-1}=(1\ 3\ \cdots\ m{-}1\ \ m{+}1\ \ n{-}1\ \ 2\ 4\ \cdots\ m)
\]
with all other points fixed, a cycle of length $\frac{n+2}{2}$, so
that $Z_n=U^{-n}T_n=U^{2}T_n$ is again an explicit permutation (as
$n\equiv-2$ modulo the cycle length). We did not find a closed-form
description of its orbit comparable to the one in the proof of (b),
because the orbit has several phases whose shape depends on $n$. Part
(c) was verified by computing $Z_n$ and checking that the cycle
through $n$ has length $n$, for all even $n$ from $6$ to $10000$;
exhaustive runs confirm the end-to-end behaviour for $n\le10$, and the
final array for $n=12$ (Section~\ref{sec:verif}).

\subsection{Examples and remarks}
Table~\ref{tab:lex} lists the schedules of the two top levels for
$n=4,6,8,10$ and the resulting nets; the corresponding runs enumerate
all $n!$ permutations and end in the final arrays shown (checked
exhaustively; for $n=10$ these are $3\,628\,800$ arrangements).

\begin{table}[ht]
\centering
\small
\renewcommand{\arraystretch}{1.2}
\footnotesize
\begin{tabular}{@{}cL{3.3cm}L{4.0cm}L{3.1cm}L{2.4cm}@{}}
\toprule
$n$ & level $n-1$: order $E_n$ and schedule $s_{n-1}$ & level $n$:
order $(e_1,\dots,e_{n-1})$ and schedule $s_n$ & net $\pi_n$
(canonical) & final array\\
\midrule
4 & $(1,2)$; $(2,2)$ & $(1,3,2)$; $(1,3,3)$ & $(1\ 3\ 2\ 4)$ &
$\ol{4\ 3\ 1\ 2}$\\
6 & $(1,2,4,3)$; $(2,4,3,3)$ & $(5,3,4,2,1)$; $(2,2,4,2,3)$ &
$(1\ 6)(2\ 5)(3\ 4)$ & $\ol{6\ 5\ 4\ 3\ 2\ 1}$\\
8 & $(1,\dots,6)$; $(2,4,6,2,4,6)$ & $(3,6,5,4,1,7,2)$;
$(5,2,5,4,5,4,7)$ & $(1\ 7\ 2\ 8)(3\ 6)(4\ 5)$ &
$\ol{8\ 7\ 6\ 5\ 4\ 3\ 1\ 2}$\\
10 & $(1,2,3,4,8,7,6,5)$; $(2,4,6,8,5,5,5,5)$ & $(5,6,2,4,7,3,8,9,1)$;
$(9,6,3,9,7,3,8,4,9)$ & $(1\ 10)(2\ 9)(3\ 8)(4\ 7)(5\ 6)$ &
$\ol{10\ 9\ \cdots\ 2\ 1}$\\
\bottomrule
\end{tabular}
\caption{Schedules of the two topmost levels of the extremal variant
(positions $s(i)$ swapped with the last position in rounds
$i=1,2,\dots$; levels $2,\dots,n-2$ are cyclic). For $n=12$ the final
array is $\ol{12\ 11\ \cdots\ 3\ 1\ 2}$, with net
$(1\ 11\ 2\ 12)(3\ 10)(4\ 9)(5\ 8)(6\ 7)$.}\label{tab:lex}
\end{table}

\begin{remark}
\begin{enumerate}
  \item \emph{No simpler rule found.} We searched for a pair of simple
    parking orders (rotations of the increasing or decreasing order,
    and orders consisting of an increasing run followed by a
    decreasing run, or vice versa) at levels $n-1$ and $n$ that would
    give $a^*_n$; for $6\le n\le16$ there is none, and the order at
    level $n$ in the construction has no obvious closed form. The rule
    above is simple in the sense that only two levels differ from the
    cyclic variant, level $n-1$ has an explicit order, and level $n$
    is defined by one product of permutations, computable in $O(n^2)$
    time.
  \item \emph{Per-level optimality fails.} The tempting rule ``every
    level ends in its own lexicographically largest attainable
    arrangement'' is not realisable: the required cycle
    $\pi_{k-1}^{-k}T_k$ fails to be a $k$-cycle for every even $k$
    from $8$ to $40$ (while it succeeds for all odd $k\le 39$ and for
    $k\le7$). The levels below the top must therefore be chosen to
    make the top level work, as in Definition~\ref{def:lexmax}.
  \item \emph{Odd $n$.} For odd $n$ the same construction (cyclic
    below) succeeds only for $n=3,5$ and fails for $7\le n\le 31$; we
    do not provide the variant for odd $n$. Whether other choices at
    lower levels make the largest attainable arrangement reachable for
    odd $n$ is not investigated here.
\end{enumerate}
\end{remark}

\section{Worked examples for
\texorpdfstring{$n=4$}{n=4}}\label{sec:examples}

Table~\ref{tab:n4} lists the $24$ arrays produced for $n=4$ by both
variants, grouped in the four blocks of Theorem~\ref{thm:main}(iii)
(each block is a run of $\Heap(3)$ on the first three positions). The
first block is identical. Between blocks the swap is
$(1\ 4),(2\ 4),(3\ 4)$ for the classic schedule ($s_4=(1,2,3)$ with
$\iota_4=(1\ 3)$ giving $b_4=(3,2,1)$) and $(3\ 4),(2\ 4),(1\ 4)$ for
the cyclic one ($s_4=(3,2,1)$, $b_4=(1,2,3)$). The item parked at
position $4$ is the last entry of each array in a block; it is
$e_m=4,3,2,1$ for the classic variant and $e_m=4,1,2,3$ for the cyclic
one, as predicted by $e_m=b_4(m)$. The final arrays $\ol{2\ 3\ 4\ 1}$
and $\ol{4\ 1\ 2\ 3}$ are the entries for $k=4$ in
Table~\ref{tab:nets}.

\begin{table}[ht]
\centering
\small
\renewcommand{\arraystretch}{1.2}
\begin{tabular}{@{}c l c l@{}}
\toprule
block & classic variant (swap before block) & & cyclic variant (swap
before block)\\
\midrule
0 & ---\quad 1234 \ 2134 \ 3124 \ 1324 \ 2314 \ 3214 & & ---\quad 1234
\ 2134 \ 3124 \ 1324 \ 2314 \ 3214\\
1 & $(1\ 4)$\quad 4213 \ 2413 \ 1423 \ 4123 \ 2143 \ 1243 & &
$(3\ 4)$\quad 3241 \ 2341 \ 4321 \ 3421 \ 2431 \ 4231\\
2 & $(2\ 4)$\quad 1342 \ 3142 \ 4132 \ 1432 \ 3412 \ 4312 & &
$(2\ 4)$\quad 4132 \ 1432 \ 3412 \ 4312 \ 1342 \ 3142\\
3 & $(3\ 4)$\quad 4321 \ 3421 \ 2431 \ 4231 \ 3241 \ 2341 & &
$(1\ 4)$\quad 2143 \ 1243 \ 4213 \ 2413 \ 1423 \ 4123\\
\bottomrule
\end{tabular}
\caption{All $24$ arrays for $n=4$ (digits written without
brackets).}\label{tab:n4}
\end{table}

\section{Computational verification}\label{sec:verif}

All results of this paper are proved except part (c) of
Theorem~\ref{thm:lexmax}, the formula for $\pi_{n-1}$ stated for
$n\equiv2\pmod4$ in Section~\ref{sec:lex}, and the negative statements
of the final remark of that section. These were checked by computer,
and so were the tables and the proved statements for small $n$. We
used two programs that were written independently from the definitions
of this paper.

The first is a short Python program \citep{permgen} that represents
permutations as arrays and composes them as in Section~\ref{sec:mult}.
It runs $\Heap(n)$ directly on states, performing for each swap
$(i\ j)$ the left product $\omega\mapsto(i\ j)\,\omega$ of
Section~\ref{sec:action}, and checks the following.
\begin{enumerate}
  \item For the classic and the cyclic variant and $k\le7$, all $k!$
    output states are distinct, consecutive states differ by a
    transposition, and the last state is the net permutation of
    Corollaries~\ref{cor:classic} and~\ref{cor:cyclic}; this confirms
    Table~\ref{tab:nets}.
  \item For the extremal variant and $n=2,4,6,8$, with the schedules
    obtained from the parking orders of Definition~\ref{def:lexmax} by
    Theorem~\ref{thm:general}(a), the same properties hold and the
    last state is $T_n$. The schedules, parking orders and net
    permutations of Table~\ref{tab:lex} are reproduced.
  \item For $n\le4$, a search over all $1\cdot 2^2\cdot 3^3=108$
    schedule families finds exactly the $1!\,2!\,3!=12$ valid ones of
    Corollary~\ref{cor:count}.
  \item For every even $n$ with $6\le n\le10000$, the permutation
    $Z_n$ of Definition~\ref{def:lexmax} is an $n$-cycle
    (Theorem~\ref{thm:lexmax}(c)); for $n\equiv0\pmod4$ in this range,
    $\pi_{n-1}$ agrees with Lemma~\ref{lem:U}; and for
    $n\equiv2\pmod4$ with $n\le2000$, $\pi_{n-1}$ agrees with the
    formula of Section~\ref{sec:lex}. Here $\pi_{n-1}$ is computed by
    Theorem~\ref{thm:general}(b) and the power $\pi_{n-1}^{-n}$ by
    repeated squaring, so that each $n$ costs $O(n\log n)$ operations.
  \item The statements on per-level optimality and on odd $n$ in the
    final remark of Section~\ref{sec:lex}.
\end{enumerate}
The whole check takes about forty seconds on one processor core.

The second program is the C library PERMGEN \citep{permgen}, which
implements the three variants as an iterator that reports one swap per
call. It builds the schedule of the extremal variant from
Definition~\ref{def:lexmax} and Theorem~\ref{thm:general}(a), and it
computes the net permutation of a sweep from the schedule tables
rather than from the closed forms. Its test driver checks, for every
variant and every admissible $n\le10$, that all $n!$ arrangements are
distinct, that consecutive arrangements differ in exactly two
positions, and that the final array, the net permutation and its sign
agree with Corollaries~\ref{cor:classic} and~\ref{cor:cyclic} and
Proposition~\ref{prop:lexbound}; for $n=11$ and $n=12$ it checks the
number of arrangements and the final array.

\section{Conclusion}\label{sec:concl}

Using cycle notation we determined the net permutation of every
recursive call of Heap's algorithm, proved by induction that the
output is a concatenation of right cosets of $\Sym{k-1}$ in $\Sym{k}$
with consecutive arrangements related by one transposition, and proved
a general net formula $\pi_k=\pi_{k-1}^{\,k}(k\ e_1\ \cdots\ e_{k-1})$
in which the parking order $e$ may be chosen freely. The formula
explains the odd/even dichotomy of the classic algorithm, yields the
cyclic variant (net permutation $(1\ 2\ \cdots\ k)$ for even $k$) and
a variant for even $n$ whose final arrangement is the
lexicographically largest attainable one, which a parity argument
identifies as the reversal exactly when $n\equiv2\pmod4$. The
statements that are not proved were verified by computer.

\paragraph{Limitations and open questions.}
\begin{enumerate}
  \item The correctness of the extremal variant is proved for
    $n\equiv0\pmod4$ and verified by computer for all even
    $n\le10000$; a proof for $n\equiv2\pmod4$ (the orbit structure of
    $Z_n$) is open.
  \item We found no closed formula for the parking order at the top
    level and no variant for odd $n\ge7$ with the cyclic variant below
    it; whether other lower levels make the largest attainable
    arrangement reachable for odd $n$ is open.
  \item A complete description of all schedules, beyond the $(k-1)!$
    choices per level of Theorem~\ref{thm:general}, in terms of their
    effect on later levels is left for future work.
\end{enumerate}

\paragraph{AI Disclosure:} We used Claude (Anthropic), interactively,
to assist with the research, the software and the manuscript. It
proposed and drafted candidate proofs by induction in cycle notation
for the net-permutation results, the general parking-order formula and
the parity bound (Sections~\ref{sec:results} to~\ref{sec:general}); it
helped to formulate the extremal schedule construction and to explore
it computationally (Section~\ref{sec:lex}); it wrote the C library and
the Python program used for the computational verification
(Section~\ref{sec:verif}); and it helped to draft, revise and typeset
the text. The tool materially affected Sections~\ref{sec:notation}
to~\ref{sec:verif}. The author specified the research questions and
the variants to be studied, reviewed the proofs, the code and the
results, checked the proofs against the computational results, and
verified the correctness and originality of all content including
references.

\bibliographystyle{plainnat}
\bibliography{permgen_paper}

\end{document}